# Multi-Stage Robust Transmission Constrained Unit Commitment: A Decomposition Framework with Implicit Decision Rules

Xuan Li [1], Qiaozhu Zhai [1*], Xiaohong Guan [1,2]

[1] Systems Engineering Institute, MOE KLINNS Lab, Xi'an Jiaotong University, Xi'an, China
[2] Center for Intelligent and Networked Systems, Department of Automation, Tsinghua University, Beijing, China
[*] qzzhai@sei.xjtu.edu.cn

**HIGHLIGHTS**

- A multi-stage robust optimization-based formulation is proposed to model the TCUC problem with uncertainty.
- Nonanticipativity of the ED decisions and multi-stage robustness of the UC decisions are guaranteed.
- A decomposition framework with implicit decision rules is proposed and the formulation scale is small.
- Subproblems in the decomposition framework have a time-decoupled structure and are easy to solve.
- A numerical example is given to clarify the difference between two-stage robustness and multi-stage robustness.
- A 2383-bus system is tested and the effectiveness and efficiency of the proposed method are verified.

**ABSTRACT**

With the integration of large-scale renewable energy sources to power systems, many optimization methods have been applied to solve the stochastic/uncertain transmission-constrained unit commitment (TCUC) problem. Among all methods, two-stage and multi-stage robust optimization-based methods are the most widely adopted ones. In the two-stage methods, nonanticipativity of economic dispatch (ED) decisions are not considered. While in multi-stage methods, explicit decision rules (for example, affine decision rules) are usually adopted to guarantee nonanticipativity of ED decisions. With explicit decision rules, the computational burden can be heavy and the optimality of the solution is affected. In this paper, a multi-stage robust TCUC formulation with implicit decision rules is proposed, as well as a decomposition framework to solve it. The solutions are proved to be multi-stage robust and nonanticipativity of ED decisions is guaranteed. Meanwhile, a computationally efficient time-decoupled solution method for the feasibility check subproblems is also proposed such that the method is suitable for large-scale TCUC problems with uncertain loads/renewable injections. Numerical tests are conducted on the IEEE 118-bus system and Polish 2383-bus system. Performances of several state-of-the-art methods are compared.

**ACKNOWLEDGMENTS**

This work is supported in part by National Key R&D Project (2016YFB0901900), National Natural Science Foundations (61773309, 61773308) of China.

# 1. Introduction

Unit commitment (UC) is a fundamental task and an important tool in power system scheduling, planning, and electricity market clearing computation (Padhy, 2004). In deregulated electricity market, security (or transmission) constrained UC (SCUC/TCUC) is solved by the independent system operator (ISO) to clear the day-ahead market and to assess the system's reliability (Chen et al., 2016). Generation corporations solve price-based UC (PBUC) to maximize their profits (Tao & Shahidehpour, 2005). UC problems are also frequently combined with optimal accommodation of wind power (Tahanan, van Ackooij, Frangioni, & Lacalandra, 2014), scheduling of the energy storage (Jiang, Wang, & Guan, 2012), generation and network expansion (El-Zonkoly, 2013), etc. The aim of a typical UC problem is to find the optimal on/off schedule and generation levels of all power generators, such that the total operating cost is minimized while satisfying many kinds of system-wide and individual unit constraints(Yu et al., 2017). The most prevalent methods for UC problem with deterministic net loads are Lagrangian relaxation (LR) based methods (Bragin, Luh, Yan, Yu, & Stern, 2015; Rong, Lahdelma, & Luh, 2008) and mixed integer programming (MIP) based methods (Morales-España, Latorre, & Ramos, 2013; Yang et al., 2017).

With the large-scale integration of uncertain renewable energy sources into power systems, stochastic optimization and uncertain programming based methods for UC problem attract more attention in recent years (Delarue & D'haeseleer, 2008; Tahanan et al., 2014). The methods for UC under uncertainty can be categorized into two-stage methods and multi-stage methods (Zheng, Wang, & Liu, 2015).

The two-stage methods formulate the UC problem as a two-stage process: the (day-ahead) UC decisions are the first-stage decisions that are made before any uncertainty is unfolded. The generation dispatch decisions are the second-stage decisions that are determined based on the complete information of the uncertainty realization during the whole scheduling horizon. Typical two-stage methods include: scenario-based methods (Schulze, Grothey, & McKinnon, 2017; Zhang, Wang, Ding, & Wang, 2018), two-stage robust optimization methods (D. Bertsimas, Litvinov, Sun, Zhao, & Zheng, 2013; Jiang, Zhang, Li, & Guan, 2014), interval optimization methods (Wang, Xia, & Kang, 2011) and chance-constrained methods (H. Y. Wu, Shahidehpour, Li, & Tian, 2014).

Different from the two-stage methods, multi-stage methods consider UC under uncertainty as a multi-stage process and the formulations are in accord with the fact that the dispatch decisions are made with only the unfolded uncertainty information, in other words, *nonanticipativity* is guaranteed. Nonanticipativity means the decisions made at the current time period shouldn't rely on the unknown future realizations (Birge & Louveaux, 2011). It has been clearly pointed out in (Ben-Tal, El Ghaoui, & Nemirovski, 2009; Lorca, Sun, Litvinov, & Zheng, 2016; Qiaozhu Zhai, Li, Lei, & Guan, 2017) that: *the two-stage and the multi-stage robustness are different*. Without consideration of nonanticipativity of the dispatch decisions, UC decisions obtained by the two-stage robust method are not truly "robust" and may be infeasible in real operation. Multi-stage robust methods (reviewed in (Dimitris Bertsimas, Brown, & Caramanis, 2011; Gabrel, Murat, & Thiele, 2014)) have been used successfully in inventory control (Ben-Tal, Boaz, & Shimrit, 2009; Dimitris Bertsimas & Thiele, 2006) and production planning(Zanjani, Ait-Kadi, & Nourelfath, 2010; Zhao, Wu, & Yuan, 2016). When multi-stage robust methods are used in solving UC under uncertainty, the computational burden is the main challenge.

The representative multi-stage methods for solving UC under uncertainty include: multi-stage stochastic programming methods (MSUC) (Takriti, Birge, & Long, 1996; L. Wu, Shahidehpour, & Li, 2007), multi-stage robust methods with explicit decision rules (MRUC-E) (Lorca & Sun, 2017; Lorca et al., 2016; Warrington, Hohl, Goulart, & Morari, 2016) and all-scenario-feasible method (ASFUC) (Qiaozhu Zhai et al., 2017).

MSUC (Takriti et al., 1996; L. Wu et al., 2007) uses the distribution information of the uncertainties to construct the scenario-trees. Since there are only a finite number of scenarios included in the scenario-trees (scenario reduction is often used to reduce the scale of the trees), the feasibility of the solution to the scenarios not included in the scenario-tree is not guaranteed.

MRUC-E (multi-stage robust methods with explicit decision rules), like, on the other hand, uses uncertainty set to represent the uncertainty, and guarantees the feasibility of the solution for all possible realizations within the set. To guarantee the nonanticipativity of the dispatch decisions, some explicit decision rules are adopted, i.e., the dispatch decisions at time $t$ are formulated as certain types of functions of the uncertainty realizations from time 1 to time $t$. And in (Lorca & Sun, 2017; Lorca et al., 2016; Warrington et al., 2016), the affine functions are used. Any specific affine decision rule is an approximation to the actual optimal decision policy, so optimality of the solution is affected and feasible region is reduced. Moreover, for tractability reasons, simplified affine functions are adopted which further reduces the feasible region. (For example, in (Lorca & Sun, 2017), dispatch decisions are assumed to be affine functions with respect to the total net load demand in the current time period only.)

It is also noted that in the framework of MRUC-E, coefficients of the explicit (affine) decision rules are obtained together with the UC decisions in the first-stage, with which the dispatch decisions will be uniquely determined when the uncertainty is revealed in the following stages. However, dispatch decisions can be re-optimized by solving economic dispatch (ED) problems in a rolling manner in each time period $t$ and the pre-obtained decision rules are solely used to determine the possible range of generation levels in period $t+1$, instead of determining the actual dispatch decisions in period $t$ directly (Lorca & Sun, 2017; Lorca et al., 2016).

ASFUC (Qiaozhu Zhai et al., 2017) uses the uncertainty set and guarantees the "all-scenario-feasibility" of the UC decisions. Different from the MRUC-E, no explicit decision rules are assumed for the dispatch decisions in ASFUC. Instead, a special set of "strong nonanticipative constraints (SNC)" are constructed with coefficients obtained together with the first-stage UC decisions. The SNC provides the dispatch decisions with bounds to guarantee the feasibility for all possible future sequentially unfolded realizations. The ASFUC has a simple single-level mixed integer linear program (MILP) structure with specifically constructed scenarios, and can be directly solved by solvers like CPLEX. But the scenarios' number is close to $2^M$ ($M$ is the number of nodes with uncertain injections), and the formulation's scale will be very large for a system with large $M$.

In this paper, we propose a computationally efficient decomposition method with implicit decision rules to solve the multi-stage robust TCUC problem. The main features of the method are summarized as follows:

1) *Multi-stage robustness*: polyhedral uncertainty set is used to formulate the possible realizations of the uncertainties. And multi-stage robustness of the UC solutions and the nonanticipativity of the dispatch decisions are guaranteed. This is very important in guaranteeing the security of operations and an example is given in section 5.1.

2) *Implicit decision rules*: decision rules for the dispatch decisions are determined implicitly: a set of implicit nonanticipative constraints (NC) with novel formulation are established and integrated into the solution procedure such that the implicit decision rules are easily determined based on a rolling ED model.

3) *Computationally efficient*: a decomposition framework is established to solve the problem. Within the framework, we propose a fast time-decoupled solution method for the feasibility check subproblem. Feasibility check is a key step in all decomposition-based framework (Lee, Liu, Mehrotra, & Shahidehpour, 2014; Ye, Wang, & Li, 2016; Zeng & Zhao, 2013) and accounts for a large portion of the computational burden. The time-decoupled method can significantly reduce the computational burden and is suitable for large-scale problems.

Comprehensive numerical tests are conducted on three cases to test the performance of the proposed approach, including a Polish 2383-bus system with 2896 transmission lines.

The remainder of the paper is organized as follows: Section 2 elaborates the formulation of the multi-stage robust TCUC problem (MRUC). Section 3 reformulates the MRUC problem with implicit decision rules (MRUC-I) by introducing the nonanticipative constraints (NC). And a rolling ED model with NC is also proposed in this section. Section 4 proposes the decomposition framework to solve the MRUC-I formulation. A fast time-decoupled method for feasibility check subproblem is also proposed in section 4. Numerical testing results are discussed in section 5 and Section 6 concludes the paper.

## 2. Multi-Stage Robust TCUC

In this section, a deterministic formulation for the transmission-constrained unit commitment (TCUC) problem is given in subsection 2.1. In subsection 2.2, we generalize the deterministic TCUC to the case with uncertain net loads, and the multi-stage robust TCUC formulation is then established.

### 2.1. Deterministic Formulation

The MILP formulation of a deterministic TCUC problem is given as follows:

$$\min_{z,p} C^S(z) + C^F(p) \quad (1)$$

$$\text{s.t.} \sum_{i \in \mathcal{N}_I} p_{i,t} - \sum_{m \in \mathcal{N}_D} d_{m,t} = 0; \forall t \in \mathcal{N}_T \quad (2)$$

$$-F_l \leq \sum_{i \in \mathcal{N}_I} \Gamma^I_{l,i} p_{i,t} - \sum_{m \in \mathcal{N}_D} \Gamma^D_{l,m} d_{m,t} \leq F_l; \forall l \in \mathcal{N}_L, \forall t \in \mathcal{N}_T \quad (3)$$

$$z_{i,t} \underline{P}_i \leq p_{i,t} \leq z_{i,t} \overline{P}_i; \forall i \in \mathcal{N}_I, \forall t \in \mathcal{N}_T \quad (4)$$

$$\Delta^-_i(z_{i,t-1}, z_{i,t}) \leq p_{i,t} - p_{i,t-1} \leq \Delta^+_i(z_{i,t-1}, z_{i,t}); \forall i \in \mathcal{N}_I, \forall t \in \mathcal{N}_T \quad (5)$$

$$z \in X, z \text{ binary} \quad (6)$$

In problem (1)-(6), $\mathcal{N}_I, \mathcal{N}_D, \mathcal{N}_L, \mathcal{N}_T$ represent the sets for thermal units, loads, transmission lines, time periods, with indices of $i, m, l, t$, and cardinality of $I, M, L, T$, respectively. $z \in \{0,1\}^{I \times T}$ and $p = [p_{i,t}]^{I \times T} \in \mathbb{R}^{I \times T}$ are matrices of UC and dispatch decisions. $d = [d_{m,t}]^{M \times T} \in \mathbb{R}^{M \times T}$ is the matrix of net loads. The objective (1) is to minimize the sum of units' start-up cost $C^S(z)$ and the fuel cost $C^F(p)$. (2) is the power balance constraint. (3) is the transmission constraint based on DC flow model, where the power flow on each transmission line $l$ in each time period $t$ must be no more than the capacity of line $l$, i.e., $F_l$. $\Gamma^I_{l,i}$ and $\Gamma^D_{l,m}$ are the power transmission distribution factors (PTDF) of unit $i$ and load $m$ to line $l$, respectively. (4) represents the generation capacity constraint. (5) is the ramping constraint, where $\Delta^+$ and $\Delta^-$ are up/down ramping limits related to units' on/off status. $X$ in (6) is the feasible region for $z$ determined by constraints on UC variables including minimum on/off time constraints, must on/off constraints, etc. More detailed formulations can be found in (Carrion & Arroyo, 2006).

### 2.2. Multi-stage Robust Formulation

Without loss of generality, the nodal net loads are taken as the uncertainty. And we consider the following uncertainty set:

$$D = \{d \in \mathbb{R}^{M \times T} : d_t \in D_t, \forall t \in \mathcal{N}_T; Hd \leq G\} \quad (7)$$

(7) is a polyhedral uncertainty set, where $D_t$ is the polyhedral set for $d_t \in \mathbb{R}^{M \times 1}$ at period $t$, $Hd \leq G$ is a general compact form of the linear temporal and spatial-correlated budget constraints with $H$ and $G$ as the coefficient matrices. More details can be found in (Ye et al., 2016) and the references therein.

In real operations of power systems, UC decisions are the first-stage decisions that should be made before any uncertainty is unfolded. Then, the dispatch decision in each time period is made without complete information of the future uncertainties, i.e., the dispatch decisions are *nonanticipative*.

The key feature of the multi-stage formulation is that the nonanticipativity is considered. In fact, $p_t$ at period $t$ must depend only on the unfolded realizations up to period $t$ (denoted as $d_{[t]}=(d_1,…,d_t)$):

$$p_t = p_t(d_{[t]}) \in \mathbb{R}^{I \times 1}, \forall t \in \mathcal{N}_T \tag{8}$$

In (8), $p_t$ is assumed to be a function of $d_{[t]}$. Though the multi-stage TCUC formulations impose more computational burden compared to the two-stage ones, it is still necessary to use the multi-stage formulation to avoid infeasible UC solutions as stated in (Ben-Tal, El Ghaoui, et al., 2009; Lorca et al., 2016; Qiaozhu Zhai et al., 2017).

The multi-stage robust TCUC (MRUC) formulation is then given as follows:

**(MRUC)**

$$\min_{z, p(\cdot)} C^S(z) + C^F\left(p(d_{[t]}^{\text{rep}})\right). \tag{9}$$

$$\text{s.t.} \sum_{i \in \mathcal{N}_I} p_{i,t}(d_{[t]}) - \sum_{m \in \mathcal{N}_D} d_{m,t} = 0; \forall t \in \mathcal{N}_T, \forall d \in D \tag{10}$$

$$-F_l \leq \sum_{i \in \mathcal{N}_I} \Gamma_{l,i}^I p_{i,t}(d_{[t]}) - \sum_{m \in \mathcal{N}_D} \Gamma_{l,m}^D d_{m,t} \leq F_l; \forall l \in \mathcal{N}_L, \forall t \in \mathcal{N}_T, \forall d \in D \tag{11}$$

$$z_{i,t} \underline{P_i} \leq p_{i,t}(d_{[t]}) \leq z_{i,t} \overline{P}_i; \forall i \in \mathcal{N}_I, \forall t \in \mathcal{N}_T, \forall d \in D \tag{12}$$

$$\Delta_i^-(z_{i,t-1}, z_{i,t}) \leq p_{i,t}(d_{[t]}) - p_{i,t-1}(d_{[t]}) \leq \Delta_i^+(z_{i,t-1}, z_{i,t}); \forall i \in \mathcal{N}_I, \forall t \in \mathcal{N}_T, \forall d \in D \tag{13}$$

$$z \in X, z \text{ binary} \tag{14}$$

The objective (9) is to minimize the total costs under the representative scenario (with superscripts "rep"). The scenario with expected net load demands is often used as the representative scenario in literature. More than one representative scenarios can also be considered with the objective replaced by the weighted sum of costs under all representative scenarios. (10)-(13) must be satisfied under all possible scenarios (i.e. $\forall d \in D$) with $p_t$ assumed to be a function of $d_{[t]}$ as (8). (14) is the same as (6).

For brevity, constraints (10)-(12) is represented in a compact form as (15), where the coefficient matrices $A, B, E, F$ can be directly obtained based on the coefficients in constraints (10)-(12). (13) is represented by (16).

$$Ap_t(d_{[t]}) + Bz + Ed_t \leq F; \forall t \in \mathcal{N}_T, \forall d \in D \tag{15}$$

$$\Delta^-(z_{t-1}, z_t) \leq p_t(d_{[t]}) - p_{t-1}(d_{[t-1]}) \leq \Delta^+(z_{t-1}, z_t); \forall t \in \mathcal{N}_T, \forall d \in D \tag{16}$$

From now on, we denote the problem (9),(14)-(16) as the MRUC, which is equivalent to the problem (9)-(14).

The UC solutions of the MRUC have the property of *multi-stage robustness* as defined below.

***Definition 1*** *(Multi-stage robustness)*: A given set of UC decisions $z$ is defined as *multi-stage robust*, if there exists a set of dispatch decisions $p_t(d_{[t]}), \forall t \in \mathcal{N}_T$, such that (15) and (16) are satisfied for any possible realization $d \in D$.

## 3. MRUC with Implicit Decision Rules

The MRUC can't be directly solved for two reasons: 1) $p_t$ is an unknown function of $d_{[t]}$ as in (8), and the optimal form of this function (decision rule) is also unknown. 2) There are an infinite number of constraints (15)-(16).

In this section, we propose an implicit decision rule to reformulate the MRUC problem, so that the function $p_t(d_{[t]})$ ($\forall t$) can be determined implicitly to address difficulty 1, while the multi-stage robustness of the solution is still guaranteed. Difficulty 2 will be addressed in section 4.

In the subsection 3.1, the methods with explicit decision rules are reviewed and analysed. Subsection 3.2 reformulates the MRUC problem with the proposed implicit decision rules. Subsection 3.3 analyses the conservatism of the implicit decision rules.

### 3.1. The Case with Explicit Decision Rules

To solve the MRUC problem, explicit decision rules (affine, to be more specific) are assumed for dispatch decisions (e.g. (Lorca & Sun, 2017; Lorca et al., 2016)) in MRUC-E methods. The constraints (15)-(16) are thus transformed into linear constraints on the coefficients in the affine functions, and these constraints must be satisfied for all possible uncertain realizations included in the uncertainty set. Afterwards, the UC decisions and the coefficients of the affine decision rules are obtained together (the first stage decisions). The obtained decision rules are not directly used to determine the real dispatch decisions, and the dispatch decisions at each period are obtained by solving "policy-enforced" economic dispatch (ED) problems (where decision rules are used to determine the possible range of dispatch decisions in the next time period).

When full affine functions are adopted, the original dispatch variables $p=[p_{i,t}]^{I \times T}$ are replaced by $I \cdot T + I \cdot M \cdot T \cdot (T+1)/2$ variables, and the problem scale is too large even for a very small system (for example, the number of variables is over 100 times of the original one when $M=10$ and $T=24$). Therefore, $p_t(d_{[t]})$ in (8) is often greatly simplified as $p_t(d_t)$ in some literature ((Lorca & Sun, 2017; Lorca et al., 2016)), i.e., the dispatch decision is assumed to be related only with the uncertainties of the current time period. Even under this assumption, the number of variables is about $M$ times of the original one.

### 3.2. Reformulation with Implicit Decision Rules

In this paper, no explicit decision rules are assumed for dispatch decisions. By introducing a set of novel nonanticipative constraints, the nonanticipativity of the dispatch decisions is naturally guaranteed, and the number of variables in problem formulation is not greatly increased. The formulation proposed in this paper is named as *multi-stage robust TCUC with implicit decision rules (MRUC-I)*.

Auxiliary variables $p^{\min}, p^{\max} \in \mathbb{R}^{I \times T}$ and related constraints are introduced in MRUC-I. They are the key to guaranteeing all features of MRUC-I.

The MRUC-I formulation is given as follows:

**(MRUC-I)**

$$\min_{z, p^{rep}, p^{\max}, p^{\min}} C^S(z) + C^F(p^{rep}) \tag{17}$$

s.t. $z \in X$, $z$ binary $\tag{18}$

$$A p_t^{rep} + B z + E d_t^{rep} \leq F ; \forall t \in \mathcal{N}_T \tag{19}$$

$$p_t^{\min} \leq p_t^{rep} \leq p_t^{\max} ; \forall t \in \mathcal{N}_T \tag{20}$$

$$p_t^{\max} - p_{t-1}^{\min} \leq \Delta^+(z_{t-1}, z_t) ; \forall t \in \mathcal{N}_T \tag{21}$$

$$p_t^{\min} - p_{t-1}^{\max} \geq \Delta^-(z_{t-1}, z_t) ; \forall t \in \mathcal{N}_T \tag{22}$$

and robust constraints:
$\forall t \in \mathcal{N}_T, \forall d \in D, \exists p_t$ such that

$$A p_t + B z + E d_t \leq F \tag{23}$$

$$p_t^{\min} \leq p_t \leq p_t^{\max} \tag{24}$$

In MRUC-I, ($z$, $p^{\text{rep}}$, $p^{\text{max}}$, $p^{\text{min}}$) are the first stage decisions, and $p^{\text{rep}}$ is the ED decisions for the representative scenario. (23)-(24) are the robust constraints that must be satisfied for any $d \in D$. Similar to the MRUC-E methods, a feasibility check subproblem is solved to check whether the robust constraints (23)-(24) are satisfied. However, unlike the MRUC-E formulation, no explicit decision rules for dispatch decisions are required.

In MRUC-I, (21)-(22) and (24) are the introduced *nonanticipative constraints (NC)*. A salient feature of MRUC-I is that the multi-stage robustness (see definition 1) of UC solution is guaranteed:

**Proposition 1**: If ($z$, $p^{\text{max}}$, $p^{\text{min}}$) is a feasible solution to MRUC-I, then $z$ is multi-stage robust.

*Proof*: based on definition 1, we need only to prove that: for any $d \in D$, 1) the feasible region of $p_t$ in period $t$ depends only on $d_{[t]}$, 2) there exists $p_t$ satisfying (15) and (16).

1) With given ($z$, $p^{\text{max}}$, $p^{\text{min}}$), the feasible set of $p_t$ determined by (23)-(24) is as follows:

$$\Omega_t(z, p_t^{\text{max}}, p_t^{\text{min}}) = \left\{ p_t : Ap_t + Bz + Ed_t \leq F; p_t^{\text{min}} \leq p_t \leq p_t^{\text{max}} \right\} \quad (25)$$

It can be seen that (25) depends only on $d_t$.

2) With the given feasible solution for MRUC-I, there must exist $p_t$ satisfying (23)-(24). With 1) being proved, (15) is now guaranteed by (23), so (15) is satisfied. Meanwhile, Based on (22) and (24), we have:

$$p_t - p_{t-1} \geq p_t^{\text{min}} - p_{t-1}^{\text{max}} \geq \Delta^-(z_{t-1}, z_t) \quad (26)$$

Similarly, based on (21) and (24), we have:

$$p_t - p_{t-1} \leq p_t^{\text{max}} - p_{t-1}^{\text{min}} \leq \Delta^+(z_{t-1}, z_t) \quad (27)$$

So (16) is satisfied.

(25)-(27) mean that any $p_t$ in $\Omega_t(z, p_t^{\text{max}}, p_t^{\text{min}})$ can be defined as the function values of $p_t(d_{[t]})$ and this implies the decision rule in MRUC is not unique. Therefore, based on the idea of the "greedy algorithm", the dispatch decisions $p_t$ can be determined based on the following ED problems.

$$\min_{p_t} C^F(p_t) \quad (28)$$

$$\text{s.t. } Ap_t + Bz + Ed_t \leq F \quad (29)$$

$$p_t^{\text{min}} \leq p_t \leq p_t^{\text{max}} \quad (30)$$

Q.E.D

Note: the optimal solution to (28)-(30) depends on the load at period $t$, or in other words, optimal dispatch decision at time period $t$ is a (vector-valued) function with respect to $d_t$ and can be written as $p_t^*(d_t)$. Since the analytical expression of this function is unknown, we call it an *implicit decision rule*.

Two salient features of MRUC-I is that: 1) the number of first-stage variables (UC decisions included) is only $4I \cdot T$, much less than that of the MRUC-E with full affine decision rules; 2) the robust constraints (23)-(24) have a time-decoupled structure such that the feasibility check subproblem can be solved very efficiently.

Now, we summarize the main features of MRUC-E with affine decision rules and the MRUC-I in the following table.

Table I Comparison of MRUC-E and MRUC-I

|  | MRUC-E | MRUC-I |
| --- | --- | --- |
| 1st stage variables | $z$, all affine coefficients | $z$, $p^{\text{rep}}$, $p^{\text{max}}$, $p^{\text{min}}$ |
| Number of 1st stage variables | $2I \cdot T + I \cdot M \cdot T \cdot (T+1)/2$ (full affine decision rule)* | $4I \cdot T$ |
| ED solution | Policy-enforced robust ED (see (Lorca et al., 2016) for more details) | Implicit decision rule (28)-(30) |

* If simplified affine decision rule (like in (Lorca et al., 2016)) is adopted, the number of 1st stage variables will be reduced but the feasible region will also be reduced.

### 3.3. Is MRUC-I Too Conservative?

Proposition 1 suggests that all feasible UC solutions to MRUC-I are also feasible to MRUC. Therefore, the main concern on MRUC-I is that whether the feasible region is *reduced* by the introduced nonanticipative constraints (NC) (21)-(22) and (24). We have several observations on this issue:

Firstly, the feasible region is probably reduced in some cases, though no numerical examples have been found. This is also the case in all methods with explicit decision rules, since any specific kind of decision rule cannot cover all possibilities and the optimal choice may be missed, especially when the coefficients in the decision rules are determined in the first stage.

Secondly, it is found that in the following case, a UC solution is feasible to MRUC is also feasible to MRUC-I, i.e., the feasible region is not reduced.

***Assumption 1***: The decision rule $p_t(d_{[t]})$ in (8) and MRUC is a function (not limited to affine functions) only with respect to $d_t$, i.e., $p_t(d_{[t]})= p_t(d_t)$. (e.g. simplified decision rules adopted in (Lorca & Sun, 2017; Lorca et al., 2016)).

Under assumption 1, let $p_{i,t}(d_t)$ be the $i$-th component (function) of the vector-valued function $p_t(d_t)$, and let

$$\begin{cases} \underline{d}_t^{(i)} = \arg\min_{d \in D} p_{i,t}(d_t), \underline{p}_{i,t} = p_{i,t}(\underline{d}_t^{(i)}) \\ \overline{d}_t^{(i)} = \arg\max_{d \in D} p_{i,t}(d_t), \overline{p}_{i,t} = p_{i,t}(\overline{d}_t^{(i)}) \end{cases} \quad (31)$$

***Assumption 2***: For every $i \in \mathcal{N}_I, t \in \mathcal{N}_T$, there exist at least two possible realizations $d^{(i,t,1)}, d^{(i,t,2)} \in D$ such that

$$\begin{cases} d_t^{(i,t,1)} = \underline{d}_t^{(i)}, d_{t-1}^{(i,t,1)} = \overline{d}_{t-1}^{(i)} \\ d_t^{(i,t,2)} = \overline{d}_t^{(i)}, d_{t-1}^{(i,t,2)} = \underline{d}_{t-1}^{(i)} \end{cases} \quad (32)$$

It should be noted that assumption 2 is not strong at all, since the two realizations impose no limitations on the net load levels in periods other than $t$ and $t$-1.

***Proposition 2***: Under assumption 1 and 2, a UC solution is feasible to MRUC is also feasible to MRUC-I. In other words, the feasible region of MRUC-I is not reduced.

*Proof:* Under assumption 1 and 2, if $z$ is feasible to MRUC, then at any specific time $t$, there must exist a set of decision rules $p_t(d_t)$ and $p_{t-1}(d_{t-1})$, such that for any $d \in D$, (15) and (16) are satisfied.

We now prove that with the same $z$, at any specific time $t$, there always exists a set of $(p_t, p_t^{\max}, p_t^{\min})$, such that for any $d \in D$, (21)-(24) are satisfied.

Let $p_t = p_t(d_t)$ ($\forall t$), then for any $d \in D$, (23) is satisfied based on (15).

Let $p_{i,t}^{\max} = \overline{p}_{i,t}$ ($\forall i,t$) and $p_{i,t}^{\min} = \underline{p}_{i,t}$ ($\forall i,t$) as defined by (31). Then, (24) must be satisfied based on (31).

$$p_{i,t}^{\max} - p_{i,t-1}^{\min} = \overline{p}_{i,t} - \underline{p}_{i,t-1} = p_{i,t}(\overline{d}_t^{(i)}) - p_{i,t-1}(\underline{d}_{t-1}^{(i)}) \quad (33)$$

Based on assumption 2 and (16), we have:

$$p_{i,t}(\overline{d}_t^{(i)}) - p_{i,t-1}(\underline{d}_{t-1}^{(i)}) \leq \Delta^+(z_{t-1}, z_t) \quad (34)$$

So (21) holds. (22) can be proved similarly. Q.E.D.

Now, it is seen that the auxiliary variables $p^{\max}$, $p^{\min}$ and the introduced NC guarantee the nonanticipativity of dispatch decisions in a simple way.

Moreover, the constraint (24) has a time-decoupled structure, which is very useful in designing the solution method. Numerical testing results in section 5 suggest that by using MRUC-I, some large-scale problems can be solved efficiently and the solution quality is no worse than that of MRUC-E with affine decision rules satisfying assumption 1.

## 4. Decomposition Framework

In section 3, we reformulate the MRUC problem into MRUC-I as (17)-(24). However, there are still an infinite number of constraints (23)-(24), and MRUC-I can't be solved directly. In this section, we propose a decomposition framework to address this issue.

In subsection 4.1 the basic algorithm for the decomposition framework is presented. The solution of the subproblem is usually quite computationally demanding, so in subsection 4.2, we propose a time-decoupled and computationally efficient method for solving the subproblem. Subsection 4.3 gives two additional techniques to accelerate the overall solution procedure.

### 4.1. Basic Algorithm

The framework has a master-sub problem structure. The master problem is solved to obtain the first-stage decisions. And the decisions are checked for multi-stage robustness by the subproblem.

The master problem (MP) is given as follows:

**(MP)**

$$\min_{z, p^{rep}, p^s, p^{min}, p^{max}} C^S(z) + C^F(p^{rep}) \tag{35}$$

$$\text{s.t. (18)-(22)} \tag{36}$$

$$Ap_t^s + Bz + Ed_t^s \leq F \,; \forall t \in \mathcal{N}_T, \forall d^s \in S \tag{37}$$

$$p_t^{min} \leq p_t^s \leq p_t^{max}; \forall t \in \mathcal{N}_T, \forall d^s \in S \tag{38}$$

In (MP), set $S$ is a subset of $D$ with a finite number of elements (scenarios). The scenarios in $S$ are iteratively added during the solution process. Each set of dispatch variables $p^s \in \mathbb{R}^{I \times T}$ corresponds to a scenario $d^s \in \mathbb{R}^{M \times T}$ in $S$.

The feasibility check subproblem (SP) is as follows:

**(SP$(z, p^{max}, p^{min})$)**

$$V^{*SP} = \max_{d \in D} \min_{p, u^1, u^2} \sum_{m,t} \left( u_{m,t}^1 + u_{m,t}^2 \right) \tag{39}$$

$$\text{s.t. } Ap_t + Bz + E\left(d_t + u_t^1 - u_t^2\right) \leq F \,; \forall t \in \mathcal{N}_T \tag{40}$$

$$p_t^{min} \leq p_t \leq p_t^{max}; \forall t \in \mathcal{N}_T \tag{41}$$

$$u^1, u^2 \geq 0 \tag{42}$$

(SP) serves as the feasibility check problem for first-stage decisions ($z, p^{max}, p^{min}$) obtained by solving (MP). $u^1, u^2 \in \mathbb{R}^{M \times T}$ are non-negative slack variables. If the optimal objective value $V^{*SP} = 0$, then $z$ is multi-stage robust. Otherwise, the optimal solution $d^*$ must be added into $S$.

It is noticed that (SP($z, p^{max}, p^{min}$)) always has feasible solutions when (MP) has feasible solutions. In fact, if ($z, p^{max}, p^{min}$) is feasible to (MP), then (37)-(38) is satisfied for the scenario $d^{rep}$, so for any $d \in D$, $u^1 - u^2 = d^{rep} - d$ can always be satisfied for some nonnegative $u^1, u^2$ and thus a feasible solution to (SP) is obtained.

Solving (SP) is quite straightforward. With duality and linearization techniques, problem (39)-(42) can be transformed into an MILP problem (e.g. (D. Bertsimas et al., 2013; Jiang et al., 2012; Lee et al., 2014; Ye et al., 2016; Zeng & Zhao, 2013) with various assumptions made for the polyhedral uncertainty set). But with a large number of nodes with uncertain load/injection, this MILP problem can be very hard to solve (Ye et al., 2016). This issue will be addressed in the next subsection.

Now, the overall solution procedure to solve MRUC-I problem (17)-(24) is as follows:

**(Algorithm 1)**

**Step 1**: Initialize set $S = \emptyset$. Set tolerance $\delta$.

**Step 2**: Solve (MP) of $S$. Obtain the optimal solution $(z^*, p^{*rep}, p^{*max}, p^{*min})$. If no feasible solutions can be found, then no feasible solutions can be found for MRUC-I, stop. Otherwise, go to Step 3.

**Step 3**: Solve (SP($z^*, p^{*max}, p^{*min}$)). Obtain the optimal solution $d^*$ and the optimal objective value $V^{*SP}$.
    If $V^{*SP} < \delta$, then go to Step 4.
    Else, $S = S \cup \{d^*\}$. Go to Step 2.

**Step 4**: End the algorithm with output $(z^*, p^{*rep}, p^{*max}, p^{*min})$

## 4.2. Fast Time-Decoupled Algorithm

The feasibility check subproblem (SP) (39)-(42) can be hard to solve, especially with a large number of nodes with uncertainty ($M$ in this paper) or with long scheduling horizon (large $T$). By dualization, a bilinear term will be introduced in the objective function, and (SP) becomes a non-convex maximization problem. Many methods have been applied in solving this problem including extreme point method (EP)(Jiang et al., 2012), binary expansion method (BE)(Ye et al., 2016), KKT conditions (Zeng & Zhao, 2013), etc. One common burden is that if ramping constraints of the dispatch decisions are considered, then (SP) has time-coupled constraints and must be solved in whole.

However, in MRUC-I, the ramping constraints are only considered in (MP) for $p^{max}, p^{min}$ (see (36) or (18)-(22)). While in (SP), once (38) is satisfied (with $p^{max}, p^{min}$ obtained from (MP)), the ramping constraints are naturally satisfied (see proof of proposition 1). With this feature, we have:

***Proposition 3***: Suppose $V^{*SP}$ is the optimal value of (SP) (39)-(42), and $V_t^{*SP}$ ($\forall t \in \mathcal{N}_T$) is the optimal value of problem (SP-$t$) (43)-(46) defined below (with the same set of ($z, p^{max}, p^{min}$)):

**(SP-$t$($z, p_t^{max}, p_t^{min}$))**

$$V_t^{*SP} = \max_{d_t \in D_t} \min_{p_t, u_t^1, u_t^2} \sum_m \left(u_{m,t}^1 + u_{m,t}^2\right) \tag{43}$$

$$\text{s.t. } Ap_t + Bz + E\left(d_t + u_t^1 - u_t^2\right) \leq F \tag{44}$$

$$p_t^{min} \leq p_t \leq p_t^{max} \tag{45}$$

$$u_t^1, u_t^2 \geq 0 \tag{46}$$

Where $p_t^{max}, p_t^{min}$ are the $t$-th column of $p^{max}, p^{min}$, respectively. Then we have:

$$V^{*SP} \leq \sum_t V_t^{*SP}. \tag{47}$$

*Proof:* The sum of the optimal objective values of all (SP-$t$) is equivalent to the optimal value of the following problem:

$$\sum_t V_t^{*SP} = \max_{d \in D_1 \times \ldots \times D_T} \min_{p, u^1, u^2} \sum_{m,t} \left(u_{m,t}^1 + u_{m,t}^2\right) \tag{48}$$

s.t. (40)-(42)

Meanwhile, based on the fact of $D \subseteq D_1 \times D_2 \times \ldots \times D_T$, we have the following inequality ((40)-(42) are the constraints):

$$\max_{d \in D} \min_{p, u^1, u^2} \sum_{m,t} \left(u_{m,t}^1 + u_{m,t}^2\right) \leq \max_{d \in D_1 \times \ldots \times D_T} \min_{p, u^1, u^2} \sum_{m,t} \left(u_{m,t}^1 + u_{m,t}^2\right)$$

Then (47) is a direct result based on the above inequality and (39), (48).      Q.E.D.

Proposition 3 means: if temporal-correlated budget constraints are considered in the uncertainty set, $\sum_t V_t^{*SP}=0$ is a sufficient condition for a given $z$ to be multi-stage robust. And it is a necessary and sufficient condition when $D=D_1 \times D_2 \times \ldots \times D_T$.

With proposition 3, we can check the multi-stage robustness of the solution of (MP) by solving $T$ independent (SP-$t$) problems with much smaller scale. These problems can be solved in a parallel framework, and the overall computational burden is significantly reduced compared to the original (SP).

With separately solved (SP-$t$), it is probable that for some time periods, $V_t^{*SP} = 0$ with given $(z, p_t^{\max}, p_t^{\min})$. For these time periods, adding scenarios and constraints may have no benefits. In this sense, we should have different scenario sets in different time periods to reduce the scale of (MP). The scenario set in time $t$ is denoted as $S_t$. With $S_t$, constraints (37)-(38) in (MP) are replaced by (49)-(50).

$$\boldsymbol{A}\boldsymbol{p}_t^s + \boldsymbol{B}\boldsymbol{z} + \boldsymbol{E}\boldsymbol{d}_t^s \leq \boldsymbol{F}; \forall \boldsymbol{d}_t^s \in S_t; \forall t \in \mathcal{N}_T \tag{49}$$

$$\boldsymbol{p}_t^{\min} \leq \boldsymbol{p}_t^s \leq \boldsymbol{p}_t^{\max}; \forall \boldsymbol{d}_t^s \in S_t; \forall t \in \mathcal{N}_T \tag{50}$$

And Algorithm 1 is modified into:

**(Algorithm 2)**

**Step 1**: Initialize set $S_t = \varnothing$ ($\forall t \in \mathcal{N}_T$). Set tolerance $\delta$.

**Step 2**: Solve (MP) of $S_t$ ($\forall t \in \mathcal{N}_T$). Obtain the optimal solution $(z^*, p^{*\max}, p^{*\min})$.

**Step 3**: for $t=1:T$

    Solve (SP-$t(z^*, p_t^{*\max}, p_t^{*\min})$), and obtain the optimal $\boldsymbol{d}_t^*$ and $V_t^{*SP}$.

      if $V_t^{*SP} \geq \delta$

        set $S_t = S_t \cup \{\boldsymbol{d}_t^*\}$.

      end # *endif*

   end  #*endfor*

   if all $S_t$ are unchanged, go to Step 4.

   else, go to Step 2.

**Step 4**: End the algorithm with output $(z^*, p^{*\max}, p^{*\min})$

### 4.3. Modifications for Feasibility Check

The (SP) problem is a feasibility check problem. Once the feasibility check is passed, the multi-stage robust solution is found. In this subsection, we try to modify $(p^{\min}, p^{\max})$ in a way such that the feasibility check is easier to pass.

*1) Adding initial scenarios*

One intuitive way is to add special scenarios (like some of the scenarios correspond to the extreme points of $D$) into the initial set $S$ (or $S_t$) before solving (MP) for the first time. The solution capable of accommodating these extreme scenarios is more likely to accommodate other scenarios.

*2) Adjustment for the nonanticipative constraints*

Another way is to adjust the coefficients of the NC before checking its feasibility: in (SP), if the given $[p^{*\min}, p^{*\max}]$ has a larger range, then the feasible region of (SP) is larger, and the optimal value of (SP), i.e. $V^{*SP}$, is more likely to be zero. For this purpose, we solve the linear programming (LP) problem (51)-(56) right after solving (MP). And then replace $(p^{*\min}, p^{*\max})$ with its optimal solution for feasibility check:

$$\max_{\boldsymbol{p}^{\max}, \boldsymbol{p}^{\min}} \sum_{i,t} \left( p_{i,t}^{\max} - p_{i,t}^{\min} \right) \tag{51}$$

$$\text{s.t. } \boldsymbol{p}_t^{\max} - \boldsymbol{p}_{t-1}^{\min} \leq \Delta^+(z_{t-1}, z_t); \forall t \tag{52}$$

$$\boldsymbol{p}_t^{\min} - \boldsymbol{p}_{t-1}^{\max} \geq \Delta^-(z_{t-1}, z_t); \forall t \tag{53}$$

$$\boldsymbol{p}_t^{\min} \leq \boldsymbol{p}_t^s \leq \boldsymbol{p}_t^{\max}; \forall s, t \tag{54}$$

$$p_{i,t}^{\min} \geq z_{i,t} \underline{p}_i; \forall i, t \tag{55}$$

$$p_{i,t}^{\max} \leq z_{i,t} \overline{p}_i; \forall i, t \tag{56}$$

Problem (51)-(56) is an LP problem for $(p^{\min}, p^{\max})$, with given UC decisions $z$ and dispatch decisions $p^s$ from the optimal solution of the (MP) problem. The objective (51) is to maximize the range of $[p^{\min}, p^{\max}]$ while satisfying nonanticipative constraints (52)-(54) and generation capacity constraints (55)-(56).

The above two methods can be used together or separately from each other and in full compatibility with time-decoupled solution procedure in the preceding subsection.

## 5. Numerical Results

At the beginning of this section, a simple numerical example modified from (Qiaozhu Zhai et al., 2017) is presented to show that the feasible solution of a two-stage robust method is infeasible under some possible scenarios in real operation, while the method proposed in this paper can avoid this infeasibility. Main numerical tests are conducted in three cases: 1) IEEE 118-bus system (Shahidehpoor, Yamin, & Li, 2002) with wind power. 2) IEEE 118-bus system with a large number of loads taken as the uncertainty. 3) The Polish 2383-bus system (Zimmerman, Murillo-Sánchez, & Thomas, 2011) with load uncertainty. All numerical tests are performed with MATLAB R2015b, YALMIP toolbox (Lofberg, 2005) and CPLEX 12.5 package on an Intel Core(TM) i7-3770 CPU @ 3.40GHz PC with 8GB RAM.

### 5.1. Infeasibility of Two-stage Robust Solution: An Example

Consider a system with net load intervals as (57) (MW), and 3 thermal units (parameters are in Table II, where the notations within are the same as (1)-(6)). Unit 1 and 2 were committed at time period 0, unit 3 wasn't. Start-up cost must be included in the total cost when unit 3 is committed during the scheduling horizon.

$$D_{t=1} = 110; D_{t=2} \in [60, 160] \tag{57}$$

Table II Parameter of the units

| I | $z_{i,0}$ | $p_{i,0}$ (MW) | $\underline{p}_i$ (MW) | $\overline{p}_i$ (MW) | $\Delta_i$ (MW/h) |
|---|---|---|---|---|---|
| 1 | 1 | 80 | 40 | 130 | 30 |
| 2 | 1 | 20 | 10 | 30 | 20 |
| 3 | 0 | 0 | 10 | 30 | 20 |

*Proposition 4.* With this system settings, **UC decision given in (58) is two-stage robust.**

$$z = [z_{i,t}]^{3\times 2} = \begin{bmatrix} 1 & 1 \\ 1 & 1 \\ 0 & 0 \end{bmatrix} \tag{58}$$

*Proof:* A given UC decision $z$ is two-stage robust, if there exist feasible ED decisions $p$ for *each and every scenario* on the extreme points of the uncertainty set (D. Bertsimas et al., 2013).

In this system, the extreme scenarios are as follows (with maximum/minimum net load on $D_{t=2}$):

$$\begin{aligned} \boldsymbol{D}^1 &= (D_{t=1}^1, D_{t=2}^1) = (110, 60) \\ \boldsymbol{D}^2 &= (D_{t=1}^2, D_{t=2}^2) = (110, 160) \end{aligned} \tag{59}$$

Under UC decision (58), a group of feasible ED solutions for the two scenarios are $\boldsymbol{p}^1$ and $\boldsymbol{p}^2$ as follows:

$$\boldsymbol{p}^1 = [p_{i,t}^1]^{3\times 2} = \begin{pmatrix} 80 & 50 \\ 30 & 10 \\ 0 & 0 \end{pmatrix} \quad \boldsymbol{p}^2 = [p_{i,t}^2]^{3\times 2} = \begin{pmatrix} 100 & 130 \\ 10 & 30 \\ 0 & 0 \end{pmatrix} \tag{60}$$

The corresponding UC solution is as (58), so UC decision (58) is two-stage robust. Q.E.D.

However, (58) is not feasible for the real operation under the following circumstance:

At time period 1, we have to make ED decision $p_{t=1}$ without knowing what $D_{t=2}$ would be. Suppose, we make $p_{t=1}$ decisions as $[p^2_{i,t=1}]^{3\times1} = (100,10,0)^T$, then it is found that when time comes to period 2 and $D_{t=2}=60$, no feasible solutions can be found for $p_{t=2}$ due to ramp limit, (in fact, with the already made decision $p_{t=1}=(100,10,0)^T$, the lowest possible total power of all three units at period 2 is (100-30)+(10-0)+0=80 MW).

This infeasibility can't be identified by the two-stage robust method, and is a direct consequence of not considering the nonanticipative constraints, which requires $p^1_{t=1} = p^2_{t=1}$. But in (60), $p^1_{t=1} \neq p^2_{t=1}$.

In fact, (58) is not a multi-stage robust UC solution, which must satisfy both two-stage robustness and nonanticipative constraints (as definition 1): Under UC decisions (58), no feasible ED solutions $p^1, p^2$ (for scenarios $D^1, D^2$, respectively) can be found such that $p^1_{t=1} = p^2_{t=1}$.

Here we directly give a feasible solution of the proposed multi-stage robust method for comparison:

**UC decision as (61) is multi-stage robust:**

$$z = [z_{i,t}]^{3\times2} = \begin{bmatrix} 1 & 1 \\ 1 & 1 \\ 1 & 1 \end{bmatrix} \tag{61}$$

Where a feasible set of $p^1, p^2$ is as follows:

$$p^1 = [p^1_{i,t}]^{3\times2} = \begin{pmatrix} 70 & 40 \\ 20 & 10 \\ 20 & 10 \end{pmatrix} \quad p^2 = [p^2_{i,t}]^{3\times2} = \begin{pmatrix} 70 & 100 \\ 20 & 30 \\ 20 & 30 \end{pmatrix} \tag{62}$$

And in (62), $p^1_{t=1} = p^2_{t=1}$, so nonanticipative constraints are satisfied. And unit 3 must be committed in both periods.

By this example, it is shown that a feasible solution of a two-stage robust method can be infeasible in real operation since nonanticipativity of ED decisions are not considered in two-stage robust method. This kind of infeasibility is avoided in the proposed multi-stage robust method.

### 5.2. Case 1: IEEE 118-Bus System with Wind Power

The test system is a modified IEEE 118-bus system. There are 54 thermal units, 179 transmission lines, 90 loads and the scheduling horizon includes 24 time periods. 4 large wind farms are located on bus 10, 25, 65 and 80, with the installed capacity of 300MW each. The total installed wind capacity is 20% of the peak load, and 37.5% of the valley load. With the added wind fields, transmission capacities of all lines are expanded by 175MW.

Three methods are tested and compared on this system: 1) multi-stage robust method with explicit affine decision rules (MRUC-E) (Lorca et al., 2016), 2) the all-scenario-feasible method (ASFUC) (Qiaozhu Zhai et al., 2017), and 3) the multi-stage robust method with implicit decision rules (MRUC-I). All three methods guarantee the *multi-stage robustness* of the UC solutions.

The uncertainty set used in this test is as (63):

$$\left\{ w \in \mathbb{R}^{K\times T} \,\middle|\, w^f - \alpha(w^f - \underline{w}) \leq w \leq w^f + \alpha(\overline{w} - w^f) \right\} \tag{63}$$

Where $w$ is the matrix of uncertain wind power output with $w^f$ as the forecasted scenario. $\underline{w}$ and $\overline{w}$ are original lower and upper bounds for $w$. $\alpha \in [0,1]$ is the adjustment parameter.

The interval for total wind power output ($\alpha = 1$) divided by total loads of each time period is shown in Fig. 1.

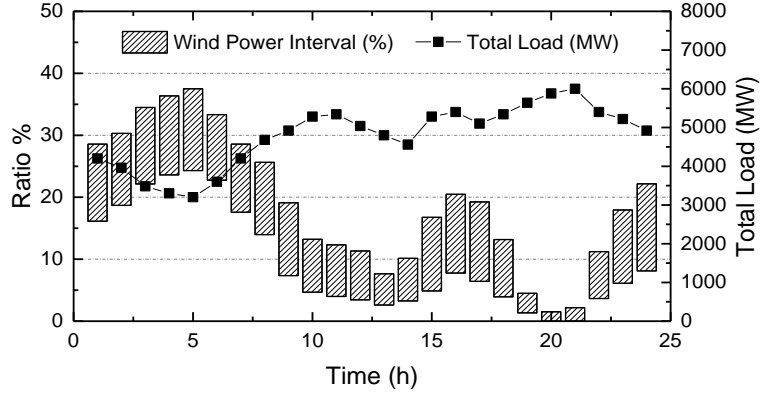

Fig. 1 Wind power output and load profile ($\alpha=1$)

Several modifications are made on these methods for comparison: 1) MRUC-E and MRUC-I are based on the same decomposition framework, where the costs under the forecasted scenario are optimized, and the sub-problem (SP) only serves as a feasibility check problem. 2) feasibility check subproblems in both MRUC-E and MRUC-I are dualized and solved using EP method as in (Jiang et al., 2012).

The affine policy used in MRUC-E is as (64), with $a$ and $b$ as the coefficients. This simplified policy achieve satisfactory performance in both optimality and computational efficiency according to the results in (Lorca et al., 2016).

$$p_{i,t} = a_{i,t}\sum_{k} w_{k,t} + b_{i,t}; \forall i,t \tag{64}$$

The optimal objective value (optimal total cost of the forecasted scenario) with different $\alpha$ are in Table III.

Table III Optimal Cost of the Forecasted Scenario ($\times 10^6$ \$)

| $\alpha$ | 1 | 0.5 | 0.3 | 0.1 | 0 |
|---|---|---|---|---|---|
| **MRUC-E** | Infeasible | Infeasible | 1.3931 | 1.3930 | 1.3722 |
| **ASFUC** | 1.3832 | 1.3726 | 1.3722 | 1.3722 | 1.3722 |
| **MRUC-I** | 1.3831 | 1.3727 | 1.3722 | 1.3722 | 1.3722 |

It is seen from Table III that: 1) no feasible solutions can be found by MRUC-E when the wind uncertainty is significant ($\alpha$=0.5, 1). And when feasible solutions can be found by MRUC-E ($\alpha$=0.3, 0.1), the optimal objective values are larger than those of ASFUC and MRUC-I. This is because the affine policy is an approximation to the full adaptive policy, and the corresponding feasible region is reduced. With the reduced feasible region, optimality of the solution is also affected. 2) as the wind uncertainty increases ($\alpha$ gets larger) the cost for the same forecasted scenario increases. This is because more units are committed to accommodating the possible wind fluctuation. Some of the units thus run in the less efficient working conditions. 3) For the case with no uncertainty ($\alpha$=0), the optimal values are the same.

The computational time and the number of scenarios of the master problem (MP) in the final iteration (in ASFUC method, only one MP needs to be solved.) are presented in Table IV.

Table IV Computational Time (s) and Final Number of Scenarios in MP

| $\alpha$ | 1 | 0.5 | 0.3 | 0.1 | 0 |
|---|---|---|---|---|---|
| **MRUC-E** | Infeasible | Infeasible | 96.8(10) | 72.8(8) | 4.5(1) |
| **ASFUC** | 1772.2(17) | 611.6(17) | 503.1(17) | 609.2(17) | 6.3(1) |
| **MRUC-I** | 92.5(5) | 45.7(4) | 72.4(5) | 29.9(3) | 6.2(1) |

It is seen from Table IV that: 1) ASFUC method is the slowest among all. Total computational time of MRUC-E and MRUC-I are close. 2) ASFUC has the most scenarios in the MP (solved only once without iterative procedure) including 16 selected vertices scenarios (SVS) and 1 base scenario (see (Qiaozhu Zhai et al., 2017) for more details). While in the iterative framework of MRUC-I and MRUC-E, the scenarios are added iteratively, and thus no new scenarios will be added when the solution of the current (MP) is identified as multi-stage robust. 3) MRUC-I has fewer iterations than MRUC-E. This is because, in the MRUC-E, optimal UC solution and also the coefficients of the simplified affine policy ($a$, $b$ in (64)) need to be checked for robustness. And since affine policy may reduce the feasible region, finding the robust affine policies may be more difficult.

Then we conduct Monte Carlo simulations on the case with $\alpha = 0.3$ to test the economic performance of the solutions. Total 1000 wind scenarios are generated with various probability distributions for a comprehensive result. For each time period in each scenario, the policy-enforced ED problem is solved for MRUC-E as in (Lorca et al., 2016), and ED with nonanticipative constraints (28)-(30) is solved for MRUC-I and ASFUC. The results are in Table V, and the optimal objective values from the UC problems are also listed for reference.

Table V Simulated Costs with $\alpha = 0.3$ ($\times 10^6$ $)

| Methods | Objective Value | Average Cost | Worst Cost |
| --- | --- | --- | --- |
| MRUC-E | 1.3931 | 1.4119 | 1.4210 |
| ASFUC | 1.3722 | 1.3715 | 1.3799 |
| MRUC-I | 1.3722 | 1.3710 | 1.3812 |

It can be seen from Table V that: 1) MRUC-E has the largest costs of all. This again proves that the affine functions are only an approximation of the full adaptive policy and optimality of the solution is sacrificed. 2) the ASFUC and MRUC-I have close simulated costs and have better economic performance than MRUC-E.

### 5.3. Case 2: IEEE 118-Bus System with Load Uncertainty

In this subsection, we conduct tests on the original decomposition algorithm (DP) in section IV.A., the time decoupled algorithm (TD) in section IV.B., and modifications for feasibility check including adding maximum/minimum scenarios into the initial set of the master problem (IS) in section IV.C.1), and the adjustments for the nonanticipative constraints (ANC) in section IV.C.2). It should be noticed that TD is a replacement for DP, while IS and ANC work based on the decomposition framework (DP or TD), and can work together or separately.

The test system is similar to case 1, except that: 1) transmission capacities are not expanded. 2) Load demands are taken as the uncertainty with the uncertainty set as (65):

$$\left\{ \boldsymbol{d} \in \mathbb{R}^{M \times T} \,\middle|\, 0.9\boldsymbol{d}^{\mathrm{f}} \leq \boldsymbol{d} \leq 1.1\boldsymbol{d}^{\mathrm{f}} \right\} \tag{65}$$

In (65), $M$ is the number of uncertain loads chosen from the 90 loads of the original system. The rest of the load fluctuation is dealt with by spinning reserves, we use the results of spinning reserve width (Qiaozhu Zhai, Tian, & Mao, 2016) for this test.

All of the above algorithms are tested with $M$=30, 50, 70, 90 to see the general performance with different problem scale. The computational time and the number of iterations are presented in Table VI and Table VII, respectively.

Table VI Computational Time of Various Algorithms (s)

|          | M=30   | M=50  | M=70   | M=90   |
|----------|--------|-------|--------|--------|
| DP       | 120.4  | 632.8 | 1200.1 | 7034.8 |
| TD       | 53.7   | 353.7 | 1132.0 | 671.7  |
| TD+IS    | 112.18 | 141.7 | 1037.8 | 3637.3 |
| TD+ANC   | **29.7** | **77.4** | **177.8** | 715.5 |
| TD+IS+ANC| 38.2   | 221.4 | 245.5  | **450.7** |

Table VII Number of Iterations of Various Algorithms

|          | M=30 | M=50 | M=70 | M=90 |
|----------|------|------|------|------|
| DP       | 3    | 5    | 5    | 7    |
| TD       | 3    | 5    | 7    | 3    |
| TD+IS    | 3    | **2** | 6    | 9    |
| TD+ANC   | **2** | **2** | **2** | 3    |
| TD+IS+ANC| **2** | 3    | **2** | **2** |

It is seen from the above results that: 1) as $M$ increases, the computational time and number of iterations increase in general. Since more nodes with uncertainty are considered, it is harder to find the multi-stage robust solution and feasibility check subproblems are harder to solve (with larger scale). 2) TD is faster than the original DP in all cases. Though in the case of $M$=70, TD uses two more iterations to find the multi-stage robust solution, the total time is still less than DP. 3) the IS technique generally reduces the number of iterations, but the computational time is not necessarily smaller since it takes more time to solve the master problem with more (extreme) scenarios. And the algorithms with IS spend more time than the ones without IS when numbers of iterations are the same. 4) the ANC technique is very efficient. Solving linear programming problem (51)-(56) itself is highly tractable (all problems can be solved within 1 second). And TD+ANC on cases with $M$=30, 50, 70 are the fastest and have the least number of iterations among all. 5) TD+IS+ANC performs better than TD+ANC on the case with the largest $M$=90, and worse on cases with smaller M.

5.4. Case 3: Polish 2383-Bus System

In this subsection, the proposed approach is tested on Polish 2383-bus system with scheduling horizon of $T$=24. The system has 179 thermal units with linear fuel costs, 2896 transmission lines with all of the transmission constraints considered in this test. The uncertainty is assumed to be on 100 buses with uncertain loads (between 80% and 120% of the original local loads).

The computational time limit is set to 8 hours. And all of DP, TD, TD+IS, TD+ANC fail to find the multi-stage robust solutions within the time limit. The main results of TD+IS+ANC are presented in Table VIII.

Table VIII Main Results of Polish 2383-Bus System

| **Obj. Value** ($\times 10^7$ \$) | 7.0665 |
|------------------|--------|
| **CPU Time** (h) | 5.6    |
| **Iterations**   | 5      |

The result proves that TD+IS+ANC is, in fact, efficient for large-scale cases. And the performance can be further improved by adopting techniques such as the elimination of the redundant transmission constraints (Q.. Zhai, Guan, Cheng, & Wu, 2010), and parallel computing for the time-decoupled subproblems.

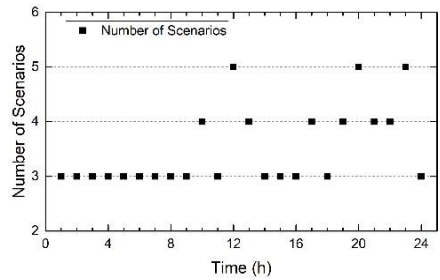

Fig. 2 Number of scenarios in each time period

The number of scenarios in each time period is presented in Fig. 2. It is seen that in the different time period, a different number of scenarios are added to the master problem in the time-decoupled framework. This prevents the method from adding unnecessary scenarios and constraints into the master problem, and hence reduces the computational burden.

## 6. Conclusions

MRUC is a more accurate formulation framework for TCUC with uncertain injections than the two-stage robust UC formulation. To solve the MRUC problem, explicit decision rules (e.g. affine functions) are usually adopted in the literature and the computational burden is heavy, unless the explicit decision rules are greatly simplified. However, it is found that the explicit decision rules can be replaced by some kind of implicit decision rules (IDR) such that much fewer decision-rule-related auxiliary variables are introduced. The IDR can guarantee the nonanticipativity of dispatch decisions and has a time-decoupled structure. With the IDR, an MRUC-I formulation is then established, and can be solved efficiently by using a decomposition method in which the time-decoupled structure of IDR plays a key role. Three cases are tested numerically, including a large-scale 2383-bus system, and the results show that the main conclusions made in this paper are valid and the proposed approach is efficient.